
\magnification=1050
\baselineskip=13pt


\def\AA{{\cal A}}  \def\BB{{\cal B}}     
\def\EE{{\cal E}}

      \def\PPP{{\bf P}} 
\def\QQQ{{\bf Q}}  \def\RRR{{\bf R}}

  \def\rightheadline{\tenrm\hfil\folio}
  \def\leftheadline{\tenrm\folio\hfil}

\def\sqr#1#2{{\vcenter{\vbox{\hrule height.#2pt
   \hbox{\vrule width.#2pt height#1pt \kern#1pt
   \vrule width.#2pt}
   \hrule height.#2pt}}}}

\def\rect#1#2#3{\raise .1ex\vbox{\hrule height.#3pt
   \hbox{\vrule width.#3pt height#2pt \kern#1pt\vrule width.#3pt}
        \hrule height.#3pt}}

\def\qed{$\hskip 5pt\rect364$} 

\def\nbar #1{{\raise7pt\hbox{\vrule height2.5pt depth-2pt width6pt}
   \kern-8pt #1}}

\def\nnbar #1{{#1\kern-7.5pt\raise6.2pt\hbox{\vrule height2.3pt depth-2pt
      width6pt}\,}}



\def\ref#1{{\rm [}{\bf #1}{\rm ]}}   
\def\nref#1#2{{\rm [}{\bf #1}{\rm ;\ #2]}}

\def\comp{{\leavevmode
     \raise.2ex\hbox{${\scriptstyle\mathchar"020E}$}}}

\def\today{\ifcase\month\or 
  January\or February\or March\or April\or  
  May\or June\or July\or August\or  
  September\or October\or November\or
  December\fi\space\number\day,\ \number\year}



\outer\def\proclaim#1{\medbreak\noindent\bf\ignorespaces
   #1\unskip.\enspace\sl\ignorespaces}
\outer\def\endproclaim{\par\ifdim\lastskip<\medskipamount\removelastskip
   \penalty 55 \fi\medskip\rm}

\font\tenBbb=msbm10

\newfam\Bbbfam \textfont\Bbbfam=\tenBbb
\font\sevenBbb=msbm7
\scriptfont8=\sevenBbb

\def\picture #1 by #2 (#3){
		\vbox to #2{
				\hrule width #1 height 0pt depth 0pt
				\vfill
				\special{picture #3}}}

\def\scaledpicture #1 by #2 (#3 scaled #4){{
		\dimen0=#1\dimen1=#2
		\divide\dimen0 by 1000\multiply\dimen0 by #4
		\divide\dimen1 by 1000\multiply\dimen1 by #4
		\picture\dimen0 by dimen1 (#3 scaled #4)}}

\outer\def\drop#1{
	\font\cap=Times at 46pt
	\setbox0=\hbox spread 2pt
	{\cap#1\hfil}
	\def\dci{\cap#1}
	\def\dc{\noindent
		\lower34.75pt\hbox
		{\hskip-\hangindent\dci}}
	\dc\vskip-34.75pt\noindent
	\hangindent
	=\wd0\hangafter=-3}

\def\pmb#1{\setbox0=\hbox{#1}%
 \kern-.025em\copy0\kern-\wd0
 \kern.05em\copy0\kern-\wd0
 \kern-.025em\raise.0433em\box0 }


\def\boxit#1{\vbox{\hrule\hbox{\vrule\kern1pt
  \vbox{\kern1pt#1\kern1pt}\kern1pt\vrule}\hrule}}

\newdimen\pxl
\pxl=.2409pt
\newbox\dlbox \setbox\dlbox=\vbox{\hrule width13\pxl height30\pxl
   depth -29\pxl \vskip-38\pxl\hrule width13\pxl height0pt depth1\pxl
   \hbox{\vrule height 30\pxl width1\pxl depth 7\pxl\kern6\pxl\vrule height
   30\pxl width1\pxl depth 7\pxl}}
\newbox\drbox \setbox\drbox=
   \vbox{\hrule width13\pxl height30\pxl depth -29\pxl
   \vskip-38\pxl\hrule width13\pxl height0pt depth 1\pxl\hbox{\kern5\pxl\vrule
   height 30\pxl width1\pxl depth 7\pxl\kern6\pxl\vrule height 30\pxl width1\pxl
   depth 7\pxl}}
\newbox\dlsbox \setbox\dlsbox=\vbox{\hrule width10\pxl height21\pxl depth
   -20\pxl \vskip-27\pxl\hrule width10\pxl height0pt depth 1\pxl\hbox{\vrule
   height 21\pxl width1\pxl depth 5\pxl\kern4\pxl\vrule height 21\pxl
   width1\pxl depth 5\pxl}}
\newbox\drsbox \setbox\drsbox=
   \vbox{\hrule width10\pxl height21\pxl depth -20\pxl
   \vskip-27\pxl\hrule width10\pxl height0pt depth
   1\pxl\hbox{\kern4\pxl\vrule height 21\pxl width1\pxl depth
   5\pxl\kern4\pxl\vrule height 21\pxl width1\pxl depth 5\pxl}}

\def\dl{\,\copy\dlbox\,} 
\def\dr{\,\copy\drbox\,} 
\def\dls{\,\copy\dlsbox\,} 
\def\drs{\,\copy\drsbox\,}
 
\def\dlb{\mathchoice{\dl}{\dl}{\dls}{\dls}}
\def\drb{\mathchoice{\dr}{\dr}{\drs}{\drs}}

\def\si#1#2 #3#4{\if#1[\dlb\else\if#1]\drb\else#1\!#1\fi\fi
   #2,#3\if#4[\dlb\else\if#4]\drb\else#4\!#4\fi\fi}
\def\grph#1#2#3{\if#1[\dlb\else\if#1]\drb\else#1\!#1\fi\fi
   #2\if#3[\dlb\else\if#3]\drb\else#3\!#3\fi\fi}

\def\pagenoslikebook{\nopagenumbers
 \headline={\ifodd\pageno\rightheadline \else\leftheadline\fi}
 \def\rightheadline{\tenrm\hfil\folio}
 \def\leftheadline{\tenrm\folio\hfil}} 

\def\m@th{\mathsurround=0pt }
\def\ialign{\everycr={}\tabskip=0pt \halign}
\def\cases#1{\left\{\,\vcenter{\normalbaselines\m@th 
\ialign{$##\hfil$&\quad##\hfil\crcr#1\crcr}}\right.}

\def\matrix#1{\null\,\vcenter{\normalbaselines\m@th
    \ialign{\hfil$##$\hfil&&\quad\hfil$##$\hfil\crcr
     \mathstrut\crcr\noalign{\kern-\baselineskip}
      #1\crcr\mathstrut\crcr\noalign{\kern-\baselineskip}}}\,}

\def\nin{\noindent}
\def\(#1){{\rm(}#1\/{\rm)}}

\def\nbar #1{{\raise7pt\hbox{\vrule height2.25pt depth-2pt width6pt}
   \kern-8pt #1}}
\def\nbara #1{{\raise4pt\hbox{\vrule height2.25pt depth-2pt width6pt}
   \kern-8pt #1}}

\def\prf{\nin{\it Proof. }}

\def\<{\langle}
\def\>{\rangle}

\def\eqdist{\,{\buildrel d\over =}\,}

\def\BOG{1}
\def\EUC{2}
\def\FER{3}
\def\LET{4}
\def\STR{5}

\centerline{\bf On a result of D.W.~Stroock}
\medskip
\centerline{P.J. Fitzsimmons}
\centerline{Department of Mathematics}
\centerline{U.C. San Diego}
\centerline{La Jolla, CA 92093--0112}
\centerline{\tt pfitzsim@ucsd.edu}
\medskip
\centerline{July 31, 2013}
\bigskip

Recently, D.W.~Stroock gave a simple probabilistic proof of L.~Schwartz' ``Borel graph theorem'', which states (in the context of Banach spaces) that if $E$ and $F$ are separable Banach spaces and $\psi:E\to F$ is a linear map  with Borel measurable graph, then  $\psi$ is continuous. In fact, Stroock obtained the continuity of $\psi$ under the weaker hypothesis that $\psi$ is $\mu$-measurable for every centered Gaussian measure $\mu$ on $E$.
  My aim here is to show that  Stroock's argument, slightly amended, proves an infinite dimensional version of the familiar fact \ref{\LET} that Lebesgue measurable solutions of Cauchy's functional equation must be continuous (and linear).

A map $\psi:E\to F$ between Banach spaces in {\it additive\/} provided  $\psi(x+y)=\psi(x)+\psi(y)$ for all $x,y\in E$. An additive $\psi$ is necessarily linear over the rationals:
$$
\psi(rx+sy) = r\psi(x)=s\psi(y),\qquad \forall r,s\in\QQQ,\forall x,y\in E,
\leqno(1)
$$

\proclaim{(2) Theorem} Let $E$ and $F$ be Banach spaces and let $\psi:E\to F$ be additive. 
If $\psi$ is $\mu$-measurable for every centered Gaussian (Radon) measure $\mu$ on $E$, then
$\psi$ is continuous (and linear).
\endproclaim

\nin (We note that in Stroock's {\it proof}, the Gaussian measures used are all Radon measures, hence our slight relaxation of his measurability assumption.)

Let us begin with a brief discussion of Gaussian measures on Banach spaces.
A probability measure $\mu$ on the Borel $\sigma$-algebra $\BB(E)$ of a Banach space $E$ is a {\it Radon\/} measure provided it is inner regular. Let $\EE:=\sigma\{ x^*: x^*\in E^*\}$ denote the cylinder $\sigma$-algebra on $E$. A Radon probability measure $\mu$ on $\BB(E)$ is uniquely determined by its restriction to $\EE$, and $\BB(E)$ is contained in the $\mu$-completion $\EE_\mu$ of $\EE$; see\nref{\BOG}{A.3.12}.

 A Radon probability measure $\mu$ on $\BB(E)$ is a centered Gaussian measure if each $x^*\in E^*$, viewed as a random variable on the probability space $(E,\EE_\mu,\mu)$, is normally distributed with mean $0$ and variance $\sigma^2(x^*)\in[0,\infty)$.
 The following characterization of Gaussian Radon measures (due to X.~Fernique) is crucial to Stroock's argument.   Let $\mu$ be a Radon probability measure on $\BB(E)$, and let $X$ and $Y$ be independent random elements of $E$ with distribution $\mu$ (defined on some probability space   
$(\Omega,\AA,\PPP)$). If $\mu$ is centered Gaussian, then for each pair of real numbers 
$(\alpha,\beta)$ with $\alpha^2+\beta^2=1$, the random vector $(\alpha X+\beta Y, \beta X-\alpha Y)$ has the same distribution as the pair $(X,Y)$, namely the product measure $\mu\otimes\mu$ (on $\EE_\mu\otimes\EE_\mu$).
(Notice that the map $(x,y)\mapsto(\alpha x+\beta y,\beta x-\alpha y)$ is $\EE_\mu\otimes\EE_\mu/\EE_\mu\otimes\EE_\mu$-measurable.) Conversely, if this equality in distribution holds for $\alpha=\beta=1/\sqrt{2}$ alone, then $\mu$ is centered Gaussian.
\medskip

\nin{\sl Proof of Theorem (2).}  The proof in \ref{\STR} needs to be supplemented at the two points where the  full linearity of $\psi$ is used: (i) in showing that the image $\psi_*\mu$ of a centered Gaussian measure $\mu$ on $E$ is a centered Gaussian measure on $F$, and (ii) in the third display on page 6 of \ref{\STR}.

Let us take up point (ii) first. The display referred to makes use of the  fact that 
$$
\<\psi(t x),y^*\>=t\< \psi(x),y^*\>,\qquad \forall t\in\RRR, x\in E, y^*\in F^*,
\leqno(3)
$$
where $F^*$ is the dual space of $F$.
To see that this partial linearity follows from our hypotheses, fix $x\in E$ and consider   the  centered Gaussian Radon measure $\mu_x$, the image of the standard normal distribution  on $\RRR$ under the mapping $\RRR\ni t\mapsto tx\in E$. The assumed $\mu_x$-measurability of $\psi$ then implies that the additive function $f(t):=\<\psi(tx),y^*\>$, $t\in\RRR$, is Lebesgue measurable. It is well known \ref{\LET} that such an $f$ is necessarily linear, and so (3) holds.

 Turning to (i), we require the following simple fact.
 
\proclaim{(4) Lemma} There is a sequence $\{(\alpha_n,\beta_n): n\ge 1\}$ of pairs of rational numbers such that $\alpha_n^2+\beta_n^2=1$ for all $n$, and $\lim_n\alpha_n=\lim_n\beta_n=1/\sqrt{2}$.
\endproclaim

\prf  We produce the required pairs by an appeal to Euclid's construction of Pythagorean triples \ref{\EUC}. Abbreviate 
$$
\kappa:=\sqrt{{\sqrt{2}-1\over\sqrt{2}+1}},
$$
choose a sequence of positive integers $m_n\in\{1,2,\ldots,n-1\}$ such that
$$
\lim_n{m_n\over n} =\kappa,
$$
and define rationals
$$
\alpha_n:={n^2-m_n^2\over n^2+m_n^2}, \qquad \beta_n:={2m_nn\over n^2+m_n^2}.
$$
Clearly $\alpha_n^2+\beta_n^2 =1$ and $\lim_n\alpha_n = (1-\kappa^2)/(1+\kappa^2) =1/\sqrt{2}$, as desired.
\qed\medskip

We now fix a centered Gaussian Radon measure $\mu$ on $\BB(E)$ and proceed to show that $\psi_*\mu$ is a centered Gaussian measure on $F$. Let $X$ and $Y$ be independent random elements of $E$ with distribution $\mu$. Let $(\alpha_n,\beta_n)$, $n\ge 1$, be as in Lemma 4. Then, using (1) for the first equality below,
$$
\eqalign{
(\alpha_n\psi(X)+\beta_n\psi(Y),\beta_n \psi(X)-\alpha_n\psi(Y))
&=(\psi(\alpha_n X+\beta_n Y),\psi(\beta_n X - \alpha_n Y))\cr
&\eqdist (\psi(X),\psi(Y)),\cr
}
$$
the $\eqdist$  indicating equality in distribution.
Sending $n\to\infty$  we obtain
$$
\left({\psi(X)+\psi(Y)\over \sqrt{2}},{\psi(X)-\psi(Y)\over\sqrt{2}}\right)\eqdist (\psi(X),\psi(Y)),
$$
so $\psi_*\mu$, the distribution of $\psi(X)$, is a centered Gaussian Radon probability measure on $\BB(F)$.
\qed
\medskip

\nin{\bf (5) Remark.} Concerning the ``universal Gaussian measurability'' hypothesis in the Theorem, we follow \ref{\STR} in noting that if $\psi:E\to F$ is additive and has a Borel measurable graph $G\subset E\times F$, and if $E$ and $F$ are  separable, then $\psi^{-1}(B) = \pi_E(G\cap (E\times B))$ is an analytic subset of $E$ (hence universally measurable) for each $B\in\BB(F)$, so Theorem 2 applies. The separability condition on $E$ is harmless, since the function $\psi$ is continuous if and only if it is sequentially continuous.
\bigskip

\centerline{\bf References}
\medskip

\baselineskip=13pt
\frenchspacing

\item{[\BOG]}
Bogachev, V.: {\it Gaussian Measures\/}, Mathematical Surveys and Monographs, {\bf 62},  American Mathematical Society, Providence,  1998.
\smallskip

\item{[\EUC]}
Euclid: {\it Elements\/}, Book X, Proposition 29, Lemma 1.
\smallskip

\item{[\FER]}
Fernique, X.:  R\'egularit\'e des trajectoires des fonctions al\'eatoires gaussiennes, In {\it  \'Ecoles d' \'et\'e de probabilit\'es de Saint-Flour IV-1974\/}, Lecture Notes in Mathematics, {\bf 480}, Springer-Verlag, Berlin, 1975, pp. 1--96.
\smallskip

\item{[\LET]}
Letac, G.: Cauchy functional equation again, {\it Amer.\ Math.\ Monthly\/} {\bf  85} (1978) 663--664.
\smallskip

\item{[\STR]}
Stroock, D.W.: On a theorem of Laurent Schwartz, {\it C.\ R.\ Acad.\ Sci.\ Paris, Ser.\ I\/}, {\bf 349} (2011) 5--6.
\smallskip

\end